\documentclass[a4paper]{article}

\usepackage{MyLaTeX}
\RequirePackage{pgf,pgffor}
\usepackage{tikz}
\tikzstyle{every node}=[circle, fill=black!0, inner sep=0pt, minimum width=4pt]

\title{The $(2,p)$-generation of sporadic simple groups}
\author{David A.~Craven, University of Birmingham}
\date{July 2015}

\begin{document}
\maketitle
\begin{abstract}
In this short note we prove that, if $p$ is an odd prime dividing the order of a sporadic simple group, then with the exception of four groups for $p=3$, all sporadic simple groups are generated by an involution and an element of order $p$.
\end{abstract}

The question of whether a finite group is generated by a particular set of elements is a natural one to ask, and in the case where the finite group is simple, there have been many papers on the subject.

All finite simple groups $G$ are generated by two elements, and given this, and the fact that `most' pairs of elements generate $G$, we can ask whether the pair of elements $(x,y)$ that generate $G$ can be given proscribed properties. One of the main such properties examined in recent years as been giving the orders of $x$ and $y$. Clearly we cannot have that both $x$ and $y$ have order $2$, so the next smallest case is where $x$ has order $2$ and $y$ has order $3$; with the exception of $\PSp_4(q)$ where $q$ is a power of $2$ or $3$ (and the Suzuki groups of course), all but finitely many simple groups are generated by an element of order $2$ and one of order $3$ (see \cite{liebeckshalev1996}, \cite{luebeckmalle1999}).

We say that a group $G$ is \emph{$(r,s)$-generated} if there exists an element $x$ of order $r$ and $y$ of order $s$ such that $G=\langle x,y\rangle$. In this article we will prove the following result.

\begin{thm}
Let $G$ be a sporadic simple group. If a prime $p>2$ divides $|G|$ then $G$ is $(2,p)$-generated, unless $p=3$ and $G$ is one of the Mathieu groups $M_{11}$, $M_{22}$ and $M_{23}$, or the McLaughlin group $McL$.
\end{thm}

We prove more than this: let $\pi_{r,s}(G)$ be the probability that if $x$ is chosen at random from $G$ to have order $x$, and $y$ is chosen at random from $G$ to have order $y$, then $G=\langle x,y\rangle$, i.e., the number of generating pairs $(x,y)$ divided by the number of all pairs $(x,y)$, with $o(x)=r$ and $o(y)=s$.

We determine exact values for $\pi_{2,p}(G)$ for $G$ a Mathieu group, $HS$, the three smallest Janko groups, and $McL$. For the other sporadic groups except for the Baby Monster and Monster, we provide lower bounds. For the Baby Monster and Monster we prove that they are generated by an involution and an element from any fixed conjugacy class of elements of order at least $4$, and also point out that they are known to be $(2,3)$-generated.

The methods we use to do this are simple: for small groups we can simply use Magma or GAP to count the number of pairs of elements that generate the group, yielding an exact figure for $\pi_{2,p}(G)$.

For larger sporadic groups we use the following fact: if $G\neq \langle x,y\rangle$, then $x$ and $y$ are contained in a maximal subgroup. Let $m_{r,s}(H)$ denote the number of pairs $(x,y)$ with $o(x)=r$ and $o(x)=s$. Thus a lower bound for $\pi_{r,s}(G)$ is $1$ minus the sum of $m_{r,s}(M)$ for all maximal subgroups $M$, i.e.,
\[ 1-\frac{\sum_{M} |G:M|m_{2,p}(M)}{m_{2,p}(G)},\]
where the sum runs over all conjugacy classes of maximal subgroups of $G$. We do this for the other sporadic groups, apart from $BM$ and $M$.

One wrinkle here is that, in Magma at least, the command \texttt{MaximalSubgroups(G)} does not work for all of the sporadic groups we wish to consider, and so we must manually construct these subgroups, mostly as subgroups of $G$ but occasionally as separate entities (for example, $Fi_{23}$ inside $Fi_{24}'$). As all we are doing is counting elements there is no harm in this. We include code to construct all maximal subgroups of these groups as an appendix, but note that the author has made no attempt to optimize this, and simply used the first technique that came into the author's head. The maximal subgroups should be constructible for any of the groups given in at most a few hours of CPU time on a reasonable machine.

For $BM$ and $M$ we use the character tables, for which there is a formula to count the triples $(a,b,c)$ with $a$ and $b$ coming from fixed conjugacy classes, $c$ being a fixed element of $G$, and $ab=c$. By judicious choice of $c$ (an element with few overgroups of $\langle c\rangle$) we can prove that $G=\langle a,b\rangle$, thus proving $(r,s)$-generation if $o(a)=r$ and $o(b)=s$.

\medskip

Neither of the methods used in the larger sporadic groups are new: using maximal subgroups was the method to prove generation in \cite{liebeckshalev1996}, and the character tables was the method used in \cite{luebeckmalle1999}.

\section{Results}

For small sporadic groups we can compute an exact figure of $\pi_{2,p}(G)$, where $\pi_{2,p}(G)$ is the probability of generating a group $G$ by an element of order $2$ and an element of order $p$.

\begin{center}\begin{tabular}{lcccccccc}
\hline Group & $3$ & $5$ & $7$ & $11$ & $13$ & $17$ & $19$& $23$
\\\hline $M_{11}$ & 0 & 0.24242 & - & 0.66667 & - & - & - & -
\\ $M_{12}$ & 0.14545 & 0.22447 & - & 0.62963 & - & - & - & -
\\ $M_{22}$ & 0 & 0.25974 & 0.55758& 0.95238 & - & - & - & -
\\ $M_{23}$ & 0 & 0.15810 & 0.33570 & 0.66666 &-&-&-& 1
\\ $M_{24}$ & 0.17475 & 0.45489 & 0.64656 & 0.75261 & - & - & - & 0.90643
\\ $HS$ & 0.13601 & 0.57087 & 0.73058 & 0.93507 & - & - & - & -
\\ $J_1$ & 0.65619 & 0.73821 & 0.98565 & 0.91729 &-&-& 0.98701 & -
\\ $J_2$ & 0.44240 & 0.63492 & 0.95309 & - & - & - & - & -
\\ $J_3$ & 0.45492 & 0.93682 & - & - & - & 0.97076 & 0.98693 & -
\\ $McL$ & 0 & 0.28239 & 0.58577 & 0.89383 & - & - & - & -
\\\hline\end{tabular}\end{center}

For groups of order larger than about a billion, exact computer calculations take too long, and so we will use maximal subgroups to get lower bounds on $\pi_{2,p}(G)$, which are more often than not close to $1$ anyway.

The next table lists the quantity
\[ 1-\frac{\sum_{M} |G:M|m_{2,p}(M)}{m_{2,p}(G)},\]
where the sum runs over all conjugacy classes of maximal subgroups of $G$.

\begin{footnotesize}\begin{center}\begin{tabular}{lcccccccccc}
\hline Group & $3$ & $5$ & $7$ & $11$ & $13$ & $17$ & $19$& $23$ & $29$ & $31$
\\\hline $He$ & 0.33107 & 0.64379 & 0.96967 & - & - & 0.97319 &-&-&-&-
\\ $Suz$ & $-0.16685$ & 0.78282 & 0.87050 & 0.99625 & 0.97554 & - & - & - & - & - 
\\ $Ru$ & 0.49808 & 0.94263 & 0.98373 & - & 0.96705 & - & - & - & 0.99976 & - 
\\ $Co_3$ & 0.13421 & 0.49498 & 0.76534 & 0.93034 & - & - & - & 0.99863 & - & - 
\\ $Co_2$ & $-2.0957$ & 0.62069 & 0.85522 & 0.96718 & - & - & - & 0.99987 & - & - 
\\ $Co_1$ & 0.68246 & 0.71187 & 0.93598 & 0.98280 & 0.99650 & - & - & 0.99718 & - & - 
\\ $HN$ & 0.74712 & 0.92483 & 0.98583 & 0.99574 & - & - & 0.99995 & - & - & -
\\ $Fi_{22}$ & $-0.055978$ & 0.61634 & 0.86736 & 0.97940 & 0.98086 & - & - & - & - & -
\\ $Fi_{23}$ & 0.35844 & 0.81137 & 0.91937 & 0.98795 & 0.97927 & 0.99979 & - & 0.99996 & -& -
\\ $Fi_{24}'$ & 0.84404 & 0.92595 & 0.99598 & 0.99598 & 0.98994 & 0.99835 &-& 0.99835 & 0.99999 &-
\\ $J_4$ & 0.88160 & 0.99760 & 0.99936 & 0.99924 & - & - & -& 0.99986 & 0.99999 & 0.99998 
\\ $ON$ & 0.48140 & 0.97141 & 0.98419 & 0.99926 & - & - & 0.99362 & - & - & 0.99967
\\ $Ly$ & 0.99611 & 0.99969 & 0.99991 & 0.99999 & - & - & - & - & - & 0.99999
\\ $Th$ & 0.88284 & 0.99818 & 0.99773 & - & 0.99977 & - & 0.99998 & - & - & 0.99998
\\\hline\end{tabular}\end{center}\end{footnotesize}
For primes larger than $31$, i.e.,  $p=37$ and $G=J_4$ and $Ly$, $p=43$ and $G=J_4$, and $p=67$ and $G=Ly$, $\pi_{2,p}(G)>0.99999$.

Note the inconvenient negative signs for $p=3$ and $G=Suz$, $Co_2$ and $Fi_{22}$. Standard generators for these in the online ATLAS are of orders $2$ and $3$, and this is simply an overcounting artefact.

\medskip

For $BM$ we have a complete list of maximal subgroups but representations of them are not available in the web ATLAS, and for $M$ a complete list of maximal subgroups is not (yet) known. We will use class multiplication coefficients for these cases, that is, the number of pairs $(x,y)$ such that $xy$ is a fixed element of the group, and $x$ comes from conjugacy class $C_1$ and $y$ comes from conjugacy class $C_2$.

For $BM$ we note that the only maximal subgroup whose order is divisible by $47$ is $47\rtimes 23$. (A quick way to see this is to note that $47$ does not divide $\Phi_n(2)$ or $\Phi_n(3)$ for any $n\leq 22$ -- so it does not divide the order of the Lie type maximal subgroups -- and it does not divide the order of smaller sporadic groups, as we see from the tables above.)

The class 47A is the 172nd class in $BM$, and class 2D is the fifth. We therefore use the GAP command
\begin{center}\texttt{ ClassMultiplicationCoefficient(CharacterTable("BM"),5,i,172)}\end{center}
for all $1\leq i\leq 184$. Except for conjugacy classes of elements of orders $1$ and $2$ this class multiplication coefficient is always positive, so $BM$ is generated by an element of 2D and an element from \emph{each} conjugacy class of elements of order at least $3$. The same holds with 2C, but not with 2A or 2B, where there are a few conjugacy classes of elements of small order that do not work.

For the Monster we use $p=59$, where the only maximal subgroup (even including those that may or may not exist) of order divisible by $59$ is $\PSL_2(59)$, which has a single conjugacy class of involutions. Using the GAP command
\begin{center}\texttt{ ClassMultiplicationCoefficient(CharacterTable("M"),n,i,152)}\end{center}
where $n=2,3$ and $1\leq i\leq 194$, we see that the class multiplication coefficient is non-zero for both 2A and 2B, and all classes except for 1A, 2A, 2B, 3A, 3B, 3C and 4A. Thus, since for all $p>3$, there exists an element $t$ from 2A, an element $u$ from 2B and an element $x$ of order $p$, such that $tx=ux$ has order $59$, this means that either $G=\langle t,x\rangle$ or $G=\langle u,x\rangle$ (depending on which class intersects $\PSL_2(59)$). Of course, the Monster is $(2,3)$-generated since it is a Hurwitz group, or because so-called standard generators of $M$ have orders $2$ and $3$.

\section{Constructing maximal subgroups}

In this section we collate example programs that construct, in Magma, the maximal subgroups for the sporadic groups where the command \texttt{MaximalSubgroups(G)} does not work. These programs have not been optimized, but work in at most a few hours of CPU time, normally a lot less.

The command \texttt{MaximalSubgroups} works for the following sporadic simple groups: $M_{11}$, $M_{12}$, $M_{22}$, $M_{23}$, $M_{24}$, $HS$, $J_1$, $J_2$, $J_3$, $McL$, $He$, $Suz$, $Ru$, $Co_3$, $Co_2$ and $Fi_{22}$.

This leaves $Co_1$, $HN$, $Fi_{23}$, $Fi_{24}'$, $J_4$, $ON$, $Ly$ and $Th$, which we will do now.

\subsection{$Co_1$}

Here we assume that $G=Co_1$ is given by standard generators $x$ in class 2B and $y$ in class 3C. The largest six of these subgroups are given in the online ATLAS, and all others were constructed by the author. This code produces a set \texttt{max} of maximal subgroups of $G$.

This code required about 10 minutes to run on the author's computer.

\begin{verbatim}
OrderG:=4157776806543360000;
n2:=OrderG/89181388800+OrderG/2012774400+OrderG/389283840;
n3:=OrderG/1345036492800+OrderG/58786560+OrderG/12597120+OrderG/544320;
n5:=OrderG/3024000+OrderG/36000;
n7:=OrderG/17640+OrderG/1176;
n11:=OrderG/66;
n13:=OrderG/156;
n23:=2*OrderG/23;

max:=[]; w1:=x; w2:=y;
w3:=w1*w2; w4:=w3*w2; w5:=w3*w4; w6:=w3*w5; w7:=w6*w3; w8:=w7*w4;
w9:=w3*w8; w5:=w4^22; w6:=w3^20; w10:=w5*w6; w11:=w4^18; w1:=w10*w11;
w3:=w8^14; w4:=w3*w3; w5:=w7*w9; w7:=w5^6; w8:=w4*w7; w2:=w8*w3;
Append(~max,sub<G|w1,w2>); w1:=x; w2:=y;
w3:=w1*w2; w4:=w3*w2; w5:=w3*w4; w6:=w3*w5; w7:=w6*w3; w8:=w6*w7;
w7:=w8*w8; w6:=w7*w7; w2:=w6*w6; w5:=w4*w4; w4:=w5^-1; w6:=w4*w2;
w2:=w6*w5; w4:=w3*w3; w3:=w4^-1; w6:=w3*w1; w1:=w6*w4;
Append(~max,sub<G|w1,w2>); w1:=x; w2:=y;
w3:=w1*w2; w4:=w3*w2; w5:=w3*w4; w6:=w3*w5; w7:=w6*w3; w8:=w6*w7;
w7:=w8*w8; w2:=w7*w7; w6:=w5*w5; w1:=w5*w6; w5:=w4*w4; w6:=w4*w5;
w7:=w5*w6; w8:=w7^-1; w6:=w8*w2; w2:=w6*w7; w4:=w3*w3; w5:=w3*w4;
w6:=w5*w5; w7:=w6^-1; w8:=w7*w1; w1:=w8*w6;
Append(~max,sub<G|w1,w2>); w1:=x; w2:=y;
w3:=w1*w2; w4:=w3*w2; w5:=w3*w4; w6:=w3*w5; w1:=w3^20; w5:=w6*w4;
w2:=w5^5; w5:=w4*w4; w6:=w4*w5; w5:=w6^-1; w3:=w5*w2; w2:=w3*w6;
Append(~max,sub<G|w1,w2>); w1:=x; w2:=y;
w3:=w1*w2; w4:=w3*w2; w5:=w3*w4; w2:=w5*w5; w1:=w2*w5; w6:=w4*w5;
w2:=w4*w6; w4:=w2*w6; w6:=w3*w3; w2:=w4*w6; w4:=w2*w5; w2:=w4*w3;
Append(~max,sub<G|w1,w2>); w1:=x; w2:=y;
w3:=w1*w2; w4:=w3*w2; w5:=w4*w3; w6:=w5*w5; w1:=w5*w6;
Append(~max,sub<G|w1,w2>);

S2:=SylowSubgroup(G,2);
S3:=SylowSubgroup(G,3);
S5:=SylowSubgroup(G,5);
S7:=SylowSubgroup(G,7);

CC:=ConjugacyClasses(S2);
z:=CC[2,3];
V:=sub<G|CC[3,3],z>;
H8:=Normalizer(G,V);
for i in [7,6,5] do
  V:=sub<G|CC[3,3],CC[4,3],CC[i,3],z>; 
  H9:=Normalizer(G,V); if(Order(H9) eq 849346560) then break i; end if;
end for;

Z3:=NormalSubgroups(S3:OrderEqual:=3^6,IsAbelian:=true);
H11:=Normalizer(G,Z3[1]`subgroup);

H13:=Normalizer(G,sub<G|y>);

Z3:=NormalSubgroups(S3:OrderEqual:=3^7);
for i in Z3 do
  H15:=Normalizer(G,i`subgroup); if(Order(H15) eq 2519424) then break i; end if;
end for;

Z3:=Subgroups(S3:OrderEqual:=9);
for i in Z3 do
  if(i`length eq 3) then
    H10:=Normalizer(G,i`subgroup); if(Order(H10) eq 235146240) then break i; end if;
  end if;
end for;

Z3:=Subgroups(S3:OrderEqual:=3);
H16:=Normalizer(G,Z3[43]`subgroup);

Z16:=Subgroups(H16:OrderEqual:=Order(Alt(5)));
for i in Z16 do
  if(IsSimple(i`subgroup)) then
    GoodA5:=i`subgroup; H12:=Normalizer(G,i`subgroup); 
    if(Order(H12) eq 72576000) then break i; end if;
  end if;
end for;

Z162:=Subgroups(GoodA5:OrderEqual:=Order(Alt(4)));
H7:=Normalizer(G,Z162[1]`subgroup);


Z16:=Subgroups(H16:OrderEqual:=Order(Alt(6)));
for i in Z16 do
  if(IsSimple(i`subgroup)) then H14:=Normalizer(G,i`subgroup); break i; end if;
end for;

Z16:=Subgroups(H16:OrderEqual:=Order(Alt(7)));
H17:=Normalizer(G,Z16[1]`subgroup);

Z5:=Subgroups(S5:OrderEqual:=5);
for i in [2..9] do
  H18:=Normalizer(G,Z5[i]`subgroup);
  if(Order(H18) eq 144000) then break i; end if;
end for;

H19:=Normalizer(G,Center(S5));

Z5:=NormalSubgroups(S5:OrderEqual:=125,IsAbelian:=true);
H20:=Normalizer(G,Z5[1]`subgroup);

H21:=Normalizer(G,S7);

Z5:=Subgroups(S5:OrderEqual:=25);
for i in [2..#Z5] do
  H22:=Normalizer(G,Z5[i]`subgroup);
  if(Order(H22) eq 3000) then break i; end if;
end for;

max cat:=[H7,H8,H9,H10,H11,H12,H13,H14,H15,H16,H17,H18,H19,H20,H21,H22];
\end{verbatim}

\subsection{$HN$}

Here we assume that $G=HN$ is given by standard generators $x$ in class 2A and $y$ in class 3B. All but three of these subgroups are given in the online ATLAS, the exceptions being a subgroup of the form $(A_6\times A_6).D_8$, a $2$-local subgroup and a $3$-local subgroup. This code produces a set \texttt{max} of maximal subgroups of $G$.

This code took about 30 minutes to run on the author's computer.

\begin{verbatim}
OrderG:=273030912000000;
n2:=OrderG/177408000+OrderG/3686400;
n3:=OrderG/544320+OrderG/29160;
n5:=OrderG/630000+OrderG/500000+2*OrderG/15000+OrderG/2500;
n7:=OrderG/420;
n11:=OrderG/22;
n19:=2*OrderG/19;

max:=[];
w1:=x; w2:=y;
w3:=w1*w2; w4:=w3*w2; w5:=w3*w4; w6:=w3*w5; w8:=w6*w5;
w9:=w3*w8; w10:=w9*w4; w9:=w10*w10; w8:=w10*w9; w2:=w9*w8;
Append(~max,sub<G|w1,w2>); w1:=x; w2:=y;
w3:=w1*w2; w4:=w3*w2; w2:=w3*w4; w6:=w3*w2; w7:=w6*w3; w5:=w6^-1;
w3:=w6*w1; w1:=w3*w5; w6:=w7*w7; w5:=w6^-1; w3:=w5*w2; w2:=w3*w6;
Append(~max,sub<G|w1,w2>); w1:=x; w2:=y;
w3:=w1*w2; w4:=w3*w2; w5:=w3*w4; w6:=w3*w5; w1:=w6^10; w5:=w4^2; w6:=w5^-1;
w4:=w6*w2; w2:=w4*w5; w4:=w3^5; w5:=w4^-1; w6:=w5*w1; w1:=w6*w4;
Append(~max,sub<G|w1,w2>); w1:=x; w2:=y;
w3:=w1*w2; w4:=w3*w2; w5:=w3*w4; w6:=w3*w5; w2:=w6^5; w5:=w4^9; w6:=w5^-1;
w4:=w6*w2; w2:=w4*w5; w4:=w3^2; w5:=w4^-1; w6:=w5*w1; w1:=w6*w4;
Append(~max,sub<G|w1,w2>); w1:=x; w2:=y;
w3:=w1*w2; w4:=w3*w2; w2:=w3*w4; w6:=w3*w2; w7:=w3^9; w8:=w6^4; w9:=w8^-1; w10:=w7*w9;
w11:=w10^-1; w12:=w11*w1; w13:=w12*w10; w7:=w6*w3; w5:=w6^-1; w3:=w6*w1; w1:=w3*w5;
w6:=w7*w7; w5:=w6^-1; w3:=w5*w2; w2:=w3*w6; w3:=w1*w2; w4:=w3*w2; w5:=w3*w4; w6:=w3*w5;
w8:=w6*w5; w9:=w3*w8; w10:=w9*w4; w5:=w10^15; w6:=w8^3; w11:=w5*w6; w5:=w3*w4; w6:=w4^10;
w1:=w3^-1*w6*w3; w6:=w4^4; w7:=w5^4; w2:=w7^-1*w6*w7; w3:=w1*w2; w4:=w3*w2; w5:=w4*w2;
w1:=w5^5; w5:=w4*w3; w6:=w3*w3; w7:=w5*w6; w8:=w6*w5; w9:=w8*w7; w8:=w6*w9; w7:=w4*w8;
w4:=w3^-1; w5:=w4*w7; w12:=w11^-1; w10:=w11*w13; w9:=w10*w12; w2:=w5*w9;
Append(~max,sub<G|w1,w2>); w1:=x; w2:=y;
w3:=w1*w2; w4:=w3*w2; w5:=w3*w4; w6:=w3*w5; w22:=w6^5; w5:=w4^9; w6:=w5^-1; w4:=w6*w22;
w22:=w4*w5; w4:=w3^2; w5:=w4^-1; w6:=w5*w1; w21:=w6*w4; w3:=w21*w22; w4:=w3*w22;
w5:=w4*w4; w4:=w5^-1; w6:=w4*w22; w31:=w6*w5; w4:=w3*w3; w32:=w3*w4; w3:=w31*w32;
w33:=w31*w32; w4:=w3*w32; w5:=w3*w4; w6:=w3*w5; w7:=w6*w3; w8:=w6*w5; w9:=w3*w8;
w10:=w9*w4; w11:=w10*w4; w12:=w3*w11; w13:=w12*w3; w14:=w13*w4; w15:=w14*w4;
w24:=w7*w15; w25:=w3^4; w3:=w1*w2; w4:=w3*w2; w5:=w3*w4; w6:=w3*w5; w7:=w6*w3;
w8:=w25*w7; w9:=w7*w25; w10:=w9^-1; w9:=w10*w8; w10:=w9*w9; w2:=w7*w10; w3:=w33*w2;
w4:=w3*w2; w5:=w3*w4; w6:=w3*w5; w7:=w6*w3; w8:=w7^3; w9:=w5^13; w10:=w8*w9; w11:=w10^-1;
w3:=w24^10; w4:=w11*w3; w3:=w4*w10; w2:=w32*w3; w4:=w31*w3; w1:=w4*w3;
Append(~max,sub<G|w1,w2>); w1:=x; w2:=y;
w3:=w1*w2; w4:=w3*w2; w2:=w3*w4; w5:=w3^-1; w6:=w5*w1;
w1:=w6*w3; w3:=w4^5; w4:=w3^-1; w5:=w4*w2; w2:=w5*w3;
Append(~max,sub<G|w1,w2>); w1:=x; w2:=y;

H:=max[1];
for i in Reverse([7..12]) do
  H:=Subgroups(H:OrderEqual:=Order(H) div i)[1]`subgroup;
end for;
Append(~max,Normalizer(G,Normalizer(G,H)));
delete H;

S2:=SylowSubgroup(G,2);
CC:=ConjugacyClasses(S2);
V:=sub<S2|CC[2,3],CC[3,3],CC[4,3]>;
Append(~max,Normalizer(G,V));

w3:=w1*w2; w4:=w3*w2; w5:=w3*w4; w6:=w3*w5; w7:=w6*w3; w8:=w7*w4; w3:=w1*w6;
w4:=w3*w8; w3:=w4^19; w4:=w3^-1; w5:=w4*w2; w2:=w5*w3; w1:=w6^10;
Append(~max,sub<G|w1,w2>); w1:=x; w2:=y;
w3:=w1*w2; w4:=w3*w2; w5:=w3*w4; w6:=w3*w5; w8:=w6*w5; w9:=w3*w8; w10:=w9*w4;
w11:=w10*w4; w12:=w3*w11; w14:=w12*w5; w15:=w14*w4; w16:=w5*w15; w17:=w16*w5;
w1:=w6^10; w3:=w8*w17; w4:=w3^4; w3:=w4^-1; w5:=w4*w1; w1:=w5*w3; w3:=w17*w8;
w4:=w3^14; w3:=w4^-1; w5:=w4*w2; w2:=w5*w3;
Append(~max,sub<G|w1,w2>); w1:=x; w2:=y;
w3:=w1*w2; w4:=w3*w2; w5:=w3*w4; w6:=w3*w5; w7:=w6*w3; w8:=w6*w5; w9:=w3*w8;
w10:=w9*w4; w11:=w10*w4; w12:=w3*w11; w14:=w12*w5; w15:=w14*w4; w16:=w5*w15;
w17:=w16*w5; w1:=w6^10; w3:=w7*w17; w4:=w3^8; w3:=w4^-1; w5:=w4*w1; w1:=w5*w3;
w3:=w17*w7; w4:=w3^9; w3:=w4^-1; w5:=w3*w2; w2:=w5*w4;
Append(~max,sub<G|w1,w2>); w1:=x; w2:=y;

S3:=SylowSubgroup(G,3);
ZS3:=NormalSubgroups(S3:OrderEqual:=81,IsElementaryAbelian:=true);
Append(~max,Normalizer(G,ZS3[1]`subgroup)); delete S3; delete ZS3;

w3:=w1*w2; w4:=w3*w2; w5:=w3*w4; w6:=w3*w5; w1:=w6^10; w5:=w4^8; w6:=w5^-1;
w4:=w6*w2; w2:=w4*w5; w4:=w3^9; w5:=w4^-1; w6:=w5*w1; w1:=w6*w4;
Append(~max,sub<G|w1,w2>);
\end{verbatim}

\subsection{$Fi_{23}$}

Here we assume that $G=Fi_{23}$ is given by standard generators $x$ in class 2B and $y$ in class 3D. All but four of these subgroups are given in the online ATLAS, the exceptions being two $3$-local subgroups and two $2$-local subgroups (the second masquerading as $\PSp_6(2)\times S_4$, so doesn't immediately look like a $2$-local). This code produces a set \texttt{max} of maximal subgroups of $G$.

This code took less than five minutes to run on the author's machine.

\begin{verbatim}
max:=[]; w1:=x; w2:=y;
w3:=w1*w2; w4:=w3*w2; w5:=w4*w4; w4:=w5*w5; w2:=w3*w3; w5:=w2*w4; w2:=w5*w3;
Append(~max,sub<G|w1,w2>); w1:=x; w2:=y;
w3:=w1*w2; w4:=w3*w2; w5:=w3^3; w6:=w4^3; w4:=w6*w5;
w1:=w4*w4; w5:=w3^8; w3:=w5*w4; w2:=w3*w5;
Append(~max,sub<G|w1,w2>); w1:=x; w2:=y;
w3:=w1*w2; w4:=w3^12; w3:=w2*w4; w2:=w3*w3;
Append(~max,sub<G|w1,w2>); w1:=x; w2:=y;
w3:=w1*w2; w4:=w3*w2; w5:=w4*w4; w2:=w5*w4; w4:=w2*w3; w5:=w3*w3; w2:=w5*w4;
Append(~max,sub<G|w1,w2>); w1:=x; w2:=y;
w3:=w1*w2; w4:=w3*w2; w5:=w3^3; w2:=w4^3; w6:=w2*w5; w1:=w6*w6; w5:=w3*w4;
w6:=w3*w5; w3:=w5*w4; w5:=w6*w3; w4:=w3*w6; w2:=w5*w4;
Append(~max,sub<G|w1,w2>); w1:=x; w2:=y;
w3:=w1*w2; w4:=w2*w1; w5:=w3^7; w3:=w5*w4; w2:=w3*w4;
Append(~max,sub<G|w1,w2>); w1:=x; w2:=y;

S3:=SylowSubgroup(G,3);
Append(~max,Normalizer(G,Center(S3)));
ZS3:=NormalSubgroups(S3:OrderEqual:=3^10);
for i in ZS3 do
  if(Order(Normalizer(G,i`subgroup)) eq 663238368) then
    Append(~max,Normalizer(G,i`subgroup)); break i;
  end if;
end for;

w3:=w1*w2; w4:=w3*w2; w5:=w3^9; w3:=w4^12; w4:=w5*w3; w2:=w4*w5;
Append(~max,sub<G|w1,w2>); w1:=x; w2:=y;
w3:=w1*w2; w4:=w3*w2; w5:=w3^3; w6:=w4^3; w2:=w6*w5; w1:=w2*w2; w5:=w3^11;
w6:=w4^13; w4:=w3*w6; w3:=w5*w6; w2:=w4*w3;
Append(~max,sub<G|w1,w2>); w1:=x; w2:=y;

S2:=SylowSubgroup(G,2);
ZS2:=NormalSubgroups(S2:OrderEqual:=2^14);
for i in ZS2 do
  if(Order(Normalizer(G,i`subgroup)) eq 247726080) then
    Append(~max,Normalizer(G,i`subgroup)); break i;
  end if; 
end for;

H:=max[9];
for i in Reverse([5..12]) do
  H:=Subgroups(H:OrderEqual:=Order(H) div i)[1]`subgroup;
end for;
Append(~max,Normalizer(G,H));
delete H;

w3:=w1*w2; w4:=w2*w1; w5:=w4^12; w2:=w3^16; w6:=w3*w5; w3:=w6*w2; w6:=w3*w5; w2:=w6*w4;
Append(~max,sub<G|w1,w2>); w1:=x; w2:=y;
w3:=w1*w2; w4:=w3*w2; w5:=w3^3; w1:=w4^10; w6:=w5*w1; w1:=w6*w3; w3:=w4^3;
w4:=w3*w5; w6:=w1^-1; w3:=w2*w1; w2:=w6*w3; w1:=w4*w4;
Append(~max,sub<G|w1,w2>);
\end{verbatim}

\subsection{$Fi_{24}'$}

Here we assume that $G=Fi_{24}'$ is given by standard generators $x$ in class 2A and $y$ in class 3E. As not all maximal subgroups will be constructed inside $G$, we require oracles for \texttt{fi23} and \texttt{he2}, the latter being the extension $He:2$. It produces \texttt{max1}, a list of maximals constructed without $G$, and \texttt{max2}, a list of maximal subgroups of $G$.

There are no maximal subgroups of $G$ given in the online ATLAS, and indeed the maximal subgroups are not even listed in the latest version there, and we have to rely on Version 2 at the time of writing.

This code took about 90 minutes to run on the author's computer.

\begin{verbatim}
OrderG:=1255205709190661721292800;
n2:=OrderG/258247006617600+OrderG/160526499840;
n3:=OrderG/44569618329600+OrderG/2424391326720
  +OrderG/14285134080+OrderG/153055008+OrderG/38211264;
n5:=OrderG/907200;
n7:=OrderG/17640+OrderG/2058;
n11:=OrderG/132;
n13:=OrderG/234;
n17:=OrderG/17;
n23:=2*OrderG/23;
n29:=2*OrderG/29;

L283:=PermutationGroup<9|\[2,1,5,6,3,4,9,8,7],\[1,3,4,2,5,7,8,6,9]>;
G1:=DirectProduct(Sym(9),Sym(5));
ZG1:=Subgroups(G1:OrderEqual:=Order(G1) div 2);
for i in ZG1 do
  if(#NormalSubgroups(i`subgroup) eq 5) then A5A92:=i`subgroup; end if;
end for;
max1:=[fi23,POmegaMinus(10,2),he2,he2,A5A92,DirectProduct(Alt(7),AGL(1,7)),
  DirectProduct(L283,Alt(6)),PGammaU(3,3),PGammaU(3,3),PGL(2,13),PGL(2,13)];
G1:=AGL(1,29); Append(~max1,Subgroups(G1:OrderEqual:=29*14)[1]`subgroup);
// max1 contains maximals not embedded in G
// These are numbers 1,4,12,13,17,18,20,21,22,23,24,25

H2:=Centralizer(G,x);
repeat a:=Random(G); until Order(a) eq 60;
b:=a^20;
H3:=Normalizer(G,sub<G|b>);

HH2:=DerivedSubgroup(H2);
HH2S2:=SylowSubgroup(HH2,2);
CCHHS2:=ConjugacyClasses(HH2S2);
for i in Reverse([5..8])
  do A:=sub<G|x,CCHHS2[i,3]>; H8:=Normalizer(G,A); 
  if(Order(H8) eq 220723937280) then break i; end if;
end for;

ZH2:=MaximalSubgroups(H2);
for i in ZH2 do
  if(Order(i`subgroup) eq 4180377600) then 
    H11:=Normalizer(G,DerivedSubgroup(DerivedSubgroup(i`subgroup))); break i;
  end if;
end for;

XXH3:=ConjugacyClasses(H3);
for i in XXH3 do
  if(Order(H3) div i[2] eq 38211264) then
    H16:=Normalizer(G,Centralizer(H3,i[3]));
    break i;
  end if;
end for;
    
S2:=SylowSubgroup(G,2);
CC2:=ConjugacyClasses(S2);
H9:=Centralizer(G,CC2[99,3]);

T1:=[CC2[i,3]:i in [11..17]];
for i in T1 do
  A:=sub<S2|Conjugates(S2,i)>;
  if(Order(A) eq 64) then
    H15:=Normalizer(G,A);
    if(Order(H15) eq 1981808640) then break i; end if;
  end if;
end for;

for i in Reverse([84..95]) do
  A:=sub<S2|Conjugates(S2,CC2[i,3])>;
  if(Order(A) eq 2^11) then H7:=Normalizer(G,A); break i; end if;
end for;

for i in [5,4,6] do
  A:=sub<G|CC2[2,3],CC2[3,3],CC2[i,3]>; H14:=Normalizer(G,A);
  if(Order(H14) eq 1981808640) then break i; end if;
end for;

S3:=SylowSubgroup(G,3);
H6:=Normalizer(G,Center(S3));
CC3:=ConjugacyClasses(S3);
A:=sub<G|CC3[2,3],CC3[4,3]>; H10:=Normalizer(G,A);

// The loop here is to remove elements that do not commute with their conjugates

T1:=[CC3[i,3]:i in [90..271]];
T2:=[];
for i in T1 do
  for j in [1..10] do g:=Random(S3);
    if((i,i^g) ne S3!1) then continue i; end if; end for;
  Append(~T2,i);
end for;
for i in Reverse(T1) do
  A:=sub<S3|Conjugates(S3,i)>;
  if(Order(A) eq 2187) then
    H5:=Normalizer(G,A);
    if(Order(H5) eq 10028164124160) then break i; end if;
  end if;
end for;

T1:=[CC3[i,3]:i in [6..13]];
for i in T1 do
  A:=sub<S3|Conjugates(S3,i)>;
  if(Order(A) eq 27) then
    H19:=Normalizer(G,A);
    if(Order(H19) eq 17907435936) then break i; end if;
  end if;
end for;

// max2 contains maximals embedded in G
// These are numbers 2,3,5,6,7,8,9,10,11,14,15,16,19

max2:=[H2,H3,H5,H6,H7,H8,H9,H10,H11,H14,H15,H16,H19];
\end{verbatim}

\subsection{$J_4$}

Here we assume that $G=J_4$ is given by standard generators $x$ in class 2A and $y$ in class 4A. All of these subgroups are given in the online ATLAS. This code produces a set \texttt{max} of maximal subgroups of $G$.

This code takes seconds to run on the author's computer.

\begin{verbatim}
OrderG:=86775571046077562880;
n2:=OrderG/21799895040+OrderG/1816657920;
n3:=OrderG/2661120;
n5:=OrderG/6720;
n7:=2*OrderG/840;
n11:=OrderG/31944+OrderG/242;
n23:=OrderG/23;
n29:=OrderG div 29;
n31:=3*OrderG div 31;
n37:=3*OrderG div 37;
n43:=3*OrderG div 43;

max:=[]; w1:=x; w2:=y;
w3:=w1*w2; w4:=w3*w2; w5:=w3*w4; w6:=w3*w5; w5:=w6*w6; w6:=w5*w5;
w5:=w6*w6; w6:=w4*w2; w7:=w6^-1; w8:=w7*w5; w2:=w8*w6; w4:=w3^-1;
w5:=w4*w1; w1:=w5*w3;
Append(~max,sub<G|w1,w2>); w1:=x; w2:=y;
w3:=w1*w2; w4:=w3*w2; w5:=w3*w4; w6:=w3*w5; w1:=w5^5; w5:=w4*w2;
w2:=w6^8; w7:=w3^9; w8:=w7^-1; w9:=w8*w1; w1:=w9*w7; w6:=w5^21;
w7:=w6^-1; w8:=w7*w2; w2:=w8*w6;
Append(~max,sub<G|w1,w2>); w1:=x; w2:=y;
w3:=w1*w2; w4:=w3*w2; w5:=w4*w2; w4:=w3^8; w6:=w4^-1; w7:=w2*w4;
w2:=w6*w7; w3:=w5^4; w4:=w3^-1; w5:=w1*w3; w1:=w4*w5;
Append(~max,sub<G|w1,w2>); w1:=x; w2:=y;
w3:=w1*w2; w4:=w3^3; w5:=w4*w2; w4:=w2*w1; w1:=w5*w2; w2:=w1^3;
w1:=w4^22; w5:=w3^13; w3:=w1*w5; w5:=w4*w4; w1:=w3^-1; w4:=w1*w2;
w1:=w2*w2; w2:=w4*w3; w4:=w5^-1; w3:=w4*w1; w1:=w3*w5;
Append(~max,sub<G|w1,w2>); w1:=x; w2:=y;
w3:=w1*w2; w4:=w3*w2; w5:=w3*w4; w6:=w3*w5; w1:=w5^5; w5:=w4*w2;
w2:=w6^8; w7:=w3^5; w8:=w7^-1; w9:=w8*w1; w1:=w9*w7; w6:=w5^4;
w7:=w6^-1; w8:=w7*w2; w2:=w8*w6;
Append(~max,sub<G|w1,w2>); w1:=x; w2:=y;
w3:=w1*w2; w4:=w3*w2; w5:=w3*w4; w6:=w3*w5; w5:=w6*w4; w2:=w6^6;
w4:=w3^-1; w6:=w1*w3; w1:=w4*w6; w4:=w3^3; w6:=w5^24; w5:=w6^-1;
w6:=w4*w5; w3:=w6^-1; w4:=w3*w2; w2:=w4*w6;
Append(~max,sub<G|w1,w2>); w1:=x; w2:=y;
w3:=w1*w2; w4:=w3*w2; w5:=w3*w4; w6:=w3*w5; w5:=w6*w4; w3:=w4^4;
w4:=w5^8; w5:=w4^-1; w6:=w5*w3; w3:=w6*w4; w4:=w1*w3; w5:=w4*w3;
w6:=w5^3; w5:=w4*w4; w3:=w6*w5; w5:=w4^6; w4:=w5*w3; w3:=w4^20;
w5:=w3*w2; w6:=w2*w3; w3:=w6^-1; w6:=w3*w5; w2:=w6^3; w3:=w4*w2;
w5:=w4*w3; w3:=w5^7; w5:=w4^-1; w6:=w5*w2; w5:=w2*w6; w2:=w4^8;
w6:=w2*w5; w5:=w6^-1; w2:=w1*w4; w4:=w5*w3; w1:=w4*w6;
Append(~max,sub<G|w1,w2>); w1:=x; w2:=y;
w3:=w1*w2; w4:=w3*w2; w5:=w3*w4; w6:=w3*w5; w7:=w2*w4; w8:=w6*w2;
w1:=w8^6; w9:=w5*w5; w2:=w5*w9; w8:=w7^9; w9:=w8^-1; w7:=w9*w1;
w1:=w7*w8; w7:=w4*w6; w8:=w7^27; w3:=w6*w4; w4:=w3^49; w5:=w8*w4;
w6:=w5^-1; w7:=w6*w2; w3:=w7*w5; w2:=w3*w3;
Append(~max,sub<G|w1,w2>); w1:=x; w2:=y;
w3:=w1*w2; w4:=w3*w2; w5:=w3*w4; w6:=w3*w5; w1:=w5^5; w5:=w4*w2;
w2:=w6^8; w7:=w3^5; w8:=w7^-1; w9:=w8*w1; w1:=w9*w7; w6:=w5^20;
w7:=w6^-1; w8:=w7*w2; w2:=w8*w6;
Append(~max,sub<G|w1,w2>); w1:=x; w2:=y;
w3:=w1*w2; w4:=w3^3; w5:=w4*w2; w4:=w2*w1; w2:=w5^4; w5:=w3^16;
w6:=w5^-1; w5:=w4^24; w7:=w6*w5; w6:=w7^-1; w5:=w6*w2; w2:=w5*w7;
w5:=w4^8; w6:=w3*w5; w7:=w3^18; w5:=w6*w7; w6:=w5^-1; w7:=w6*w1;
w1:=w7*w5;
Append(~max,sub<G|w1,w2>); w1:=x; w2:=y;
w3:=w1*w2; w4:=w3^3; w5:=w4*w2; w4:=w5*w2; w1:=w4^6; w5:=w3^8;
w3:=w4*w5; w4:=w5^-1; w5:=w4*w3; w4:=w5^6; w3:=w4*w2; w6:=w2*w4;
w4:=w6^-1; w6:=w4*w3; w4:=w6^7; w3:=w2*w4; w2:=w3^-1; w4:=w5*w2;
w6:=w4*w4; w2:=w3^5; w4:=w2*w6; w6:=w4^-1; w2:=w6*w5; w5:=w2*w4;
w2:=w3*w5; w4:=w2*w2; w6:=w3*w3; w3:=w5^3; w2:=w3*w4; w7:=w4*w3;
w4:=w7^-1; w7:=w4*w2; w4:=w7*w7; w2:=w4*w6; w7:=w3*w6; w4:=w6*w3;
w3:=w4^-1; w6:=w3*w7; w4:=w2*w6; w3:=w4*w4; w2:=w5*w4; w6:=w2*w5;
w7:=w6*w2; w2:=w7*w3; w5:=w2^5; w6:=w3*w5; w7:=w3*w2; w4:=w6*w6;
w6:=w4*w7; w7:=w6^-1; w5:=w2*w6; w2:=w7*w5;
Append(~max,sub<G|w1,w2>); w1:=x; w2:=y;
w3:=w1*w2; w4:=w3*w2; w5:=w3*w4; w1:=w4*w2; w2:=w5*w5; w5:=w2*w1;
w2:=w5^3; w5:=w3*w3; w3:=w5*w1; w1:=w3^6; w3:=w4*w4; w4:=w5*w3;
w3:=w4^26; w6:=w5*w3; w3:=w6^-1; w7:=w1*w6; w1:=w3*w7; w3:=w4^33;
w6:=w5^11; w7:=w3*w6; w6:=w7^-1; w4:=w2*w7; w2:=w6*w4;
Append(~max,sub<G|w1,w2>); w1:=x; w2:=y;
w3:=w1*w2; w4:=w3*w2; w5:=w3*w4; w6:=w3*w5; w7:=w6^8; w6:=w5^5;
w8:=w3^8; w9:=w8^-1; w10:=w9*w6; w6:=w10*w8; w8:=w2*w1; w9:=w8^19;
w8:=w9^-1; w10:=w8*w7; w7:=w10*w9; w8:=w6*w7; w9:=w8^3; w10:=w1*w9;
w11:=w10^8; w3:=w2*w9; w4:=w9*w2; w5:=w3^-1; w3:=w5*w4; w4:=w3^18;
w12:=w2*w4; w1:=w11*w12; w3:=w1*w8; w4:=w3*w8; w5:=w3*w4;
w16:=w3*w5; w18:=w16*w5; w9:=w3*w18; w10:=w11*w9; w20:=w10^7;
w13:=w5^5; w14:=w8*w13; w23:=w14^3; w4:=w3^-1; w13:=w1^4; w14:=w13*w4;
w13:=w14^-1; w15:=w13*w20; w21:=w15*w14; w13:=w1^8; w14:=w4^7;
w15:=w13*w14; w14:=w15^-1; w13:=w14*w20; w22:=w13*w15; w13:=w1^10;
w14:=w4^6;w15 :=w13*w14;w14:=w15^-1; w13:=w14*w20; w24:=w13*w15; w20:=w21^-1;
w19:=w20*w8; w18:=w19*w21; w19:=w18*w8; w21:=w19*w20; w20:=w22^-1;
w19:=w20*w8; w18:=w19*w22; w19:=w18*w8; w22:=w19*w20; w20:=w24^-1;
w19:=w20*w8; w18:=w19*w24; w19:=w18*w8; w24:=w19*w20; w1:=w21*w23;
w2:=w22*w24; w3:=w1*w2; w4:=w3*w2; w5:=w4*w4; w2:=w1^-1; w3:=w5*w2;
w4:=w3^-1; w5:=w4*w21; w2:=w5*w3; w9:=w8*w7; w1:=w9*w7;
Append(~max,sub<G|w1,w2>);
\end{verbatim}

\subsection{$ON$}

Here we assume that $G=ON$ is given by standard generators $x$ in class 2A and $y$ in class 4A. All of these subgroups are given in the online ATLAS. This code produces a set \texttt{max} of maximal subgroups of $G$.

This code takes seconds to run on the author's computer.

\begin{verbatim}
OrderG:=460815505920;
n2:=OrderG/161280;
n3:=OrderG/3240;
n5:=OrderG/180;
n7:=OrderG/1372+OrderG/49;
n11:=OrderG/11;
n19:=3*OrderG/19;
n31:=2*OrderG/31;

max:=[]; w1:=x; w2:=y;
w3:=w1*w2; w4:=w3*w2; w5:=w4*w4; w4:=w5^-1; w6:=w4*w2;
w2:=w6*w5; w4:=w3^-1; w5:=w4*w1; w1:=w5*w3;
Append(~max,sub<G|w1,w2>); w1:=x; w2:=y;
w3:=w1*w2; w4:=w3*w2; w5:=w4*w4; w6:=w5^-1; w3:=w2^-1;
w4:=w2*w1; w1:=w4*w3; w4:=w6*w3; w2:=w4*w5;
Append(~max,sub<G|w1,w2>); w1:=x; w2:=y;
w3:=w1*w2; w4:=w3*w2; w5:=w3*w4; w6:=w3*w5; w8:=w6*w5; w9:=w3*w8;
w8:=w9*w9; w2:=w8*w8; w6:=w5^6; w7:=w6^-1; w8:=w7*w2; w2:=w8*w6;
w5:=w4^7; w6:=w5^-1; w7:=w6*w1; w1:=w7*w5;
Append(~max,sub<G|w1,w2>); w1:=x; w2:=y;
w3:=w1*w2; w4:=w3*w2; w5:=w3*w4; w6:=w5^10; w7:=w6*w2;
w8:=w2*w6; w9:=w8^-1; w6:=w9*w7; w1:=w6^14; w2:=w3*w4;
Append(~max,sub<G|w1,w2>); w1:=x; w2:=y;
w3:=w1*w2; w4:=w3*w2; w5:=w3*w4; w6:=w3*w5; w7:=w6*w3; w8:=w7*w4; w9:=w7*w8; w2:=w9*w9;
w1:=w9*w9; w7:=w6*w4; w8:=w7^19; w9:=w8^-1; w10:=w9*w1; w1:=w10*w8; w7:=w4*w6; w8:=w7^17;
w6:=w5*w4; w7:=w6*w3; w9:=w7^8; w7:=w8*w9; w6:=w7^-1; w8:=w6*w2; w2:=w8*w7;
Append(~max,sub<G|w1,w2>); w1:=x; w2:=y;
w3:=w1*w2; w4:=w3*w2; w5:=w3*w4; w6:=w3*w5; w7:=w6*w3; w8:=w7*w4; w9:=w7*w8; w2:=w9*w9;
w1:=w9*w9; w7:=w6*w4; w8:=w7^19; w9:=w8^-1; w10:=w9*w1; w1:=w10*w8; w8:=w4*w6;
w6:=w5*w4; w7:=w6*w3; w9:=w7^12; w7:=w8*w9; w6:=w7^-1; w8:=w6*w2; w2:=w8*w7;
Append(~max,sub<G|w1,w2>); w1:=x; w2:=y;
w3:=w1*w2; w4:=w3*w2; w5:=w3*w4; w4:=w5*w2; w5:=w4*w4; w2:=w5*w5;
w4:=w3*w3; w5:=w3*w4; w4:=w5^-1; w3:=w4*w1; w1:=w3*w5;
Append(~max,sub<G|w1,w2>); w1:=x; w2:=y;
w3:=w1*w2; w4:=w3*w2; w5:=w4*w2; w4:=w5*w3; w3:=w4*w4; w2:=w3*w3;
w4:=w5*w5; w3:=w5*w4; w4:=w3^-1; w5:=w4*w1; w1:=w5*w3;
Append(~max,sub<G|w1,w2>); w1:=x; w2:=y;
w3:=w1*w2; w4:=w3*w2; w5:=w3*w4; w6:=w3*w5; w8:=w6*w5; w9:=w3*w8; w8:=w9*w9;
w2:=w8*w8; w5:=w4*w4; w4:=w5*w5; w5:=w4^-1; w6:=w5*w2; w2:=w6*w4;
w4:=w3*w3; w3:=w4*w4; w4:=w3^-1; w5:=w4*w1; w1:=w5*w3;
Append(~max,sub<G|w1,w2>); w1:=x; w2:=y;
w3:=w1*w2; w4:=w3*w2; w5:=w3*w4; w6:=w3*w5; w7:=w6*w3; w8:=w7*w4; w9:=w7*w8; w2:=w9*w9;
w7:=w6*w4; w8:=w7^4; w9:=w8^-1; w10:=w9*w1; w1:=w10*w8; w7:=w4*w6; w8:=w7^19;
w6:=w5*w4; w7:=w6*w3; w9:=w7^2; w7:=w8*w9; w6:=w7^-1; w8:=w6*w2; w2:=w8*w7;
Append(~max,sub<G|w1,w2>); w1:=x; w2:=y;
w3:=w1*w2; w4:=w3*w2; w5:=w3*w4; w6:=w3*w5; w7:=w6*w3; w8:=w7*w4; w9:=w7*w8;
w2:=w9*w9; w7:=w6*w4; w8:=w7^14; w9:=w8^-1; w10:=w9*w1; w1:=w10*w8; w7:=w4*w6;
w8:=w7^2; w6:=w5*w4; w9:=w6*w3; w7:=w8*w9; w6:=w7^-1; w8:=w6*w2; w2:=w8*w7;
Append(~max,sub<G|w1,w2>); w1:=x; w2:=y;
w3:=w1*w2; w4:=w3*w2; w5:=w3*w4; w6:=w3*w5; w9:=w2*w5; w2:=w9*w9; w8:=w6*w4;
w9:=w8^-1; w10:=w9*w1; w1:=w10*w8; w7:=w4*w6; w8:=w7^10; w6:=w5*w4; w7:=w6*w3;
w9:=w7^3; w7:=w8*w9; w6:=w7^-1; w8:=w6*w2; w2:=w8*w7;
Append(~max,sub<G|w1,w2>); w1:=x; w2:=y;
w3:=w1*w2; w4:=w3*w2; w5:=w3*w4; w6:=w3*w5; w9:=w2*w5; w2:=w9*w9; w7:=w6*w4;
w8:=w7^13; w9:=w8^-1; w10:=w9*w1; w1:=w10*w8; w7:=w4*w6; w8:=w7^27; w6:=w5*w4;
w7:=w6*w3; w9:=w7^3; w7:=w8*w9; w6:=w7^-1; w8:=w6*w2; w2:=w8*w7;
Append(~max,sub<G|w1,w2>);
\end{verbatim}

\subsection{$Ly$}

Here we assume that $G=Ly$ is given by standard generators $x$ in class 2A and $y$ in class 5A. All of these subgroups are given in the online ATLAS. This code produces a set \texttt{max} of maximal subgroups of $G$.

This code takes seconds to run on the author's computer.

\begin{verbatim}
OrderG:=51765179004000000
n2:=OrderG/39916800;
n3:=OrderG/2694384000+OrderG/174960;
n5:=OrderG/2250000+OrderG/3750;
n7:=OrderG/168;
n11:=OrderG/33;
n31:=5*OrderG/31;
n37:=2*OrderG/37;
n67:=3*OrderG/67;

max:=[]; w1:=x; w2:=y;
w3:=w1*w2; w4:=w3*w2; w5:=w3*w4; w9:=w5*w2; w6:=w3*w5; w7:=w6*w3;
w8:=w7*w7; w2:=w7*w8; w3:=w9^7; w4:=w3^-1; w5:=w3*w1; w1:=w5*w4;
w8:=w6^25; w7:=w8^-1; w3:=w7*w2; w2:=w3*w8;
Append(~max,sub<G|w1,w2>); w1:=x; w2:=y;
w3:=w1*w2; w4:=w3*w2; w5:=w3*w4; w9:=w2*w5; w6:=w3*w5; w7:=w6*w3;
w8:=w7*w7; w2:=w7*w8; w3:=w6^15; w4:=w3^-1; w5:=w3*w1; w1:=w5*w4;
w8:=w9^12; w7:=w8^-1; w3:=w7*w2; w2:=w3*w8;
Append(~max,sub<G|w1,w2>); w1:=x; w2:=y;
w3:=w1*w2; w4:=w3*w2; w5:=w3*w4; w6:=w3*w5; w7:=w6*w3; w5:=w3^-1; w9:=w5*w1;
w1:=w9*w3; w2:=w7^3; w6:=w4^12; w5:=w6^-1; w3:=w5*w2; w2:=w3*w6;
Append(~max,sub<G|w1,w2>); w1:=x; w2:=y;
w3:=w1*w2; w4:=w3*w2; w5:=w3*w4; w6:=w3*w5; w8:=w6*w5; w9:=w3*w8;
w10:=w9*w4;w4:=w3^7; w5:=w4*w10; w6:=w10*w4; w7:=w6^-1; w6:=w7*w5;
w7:=w6^10; w2:=w10*w7; w4:=w1*w1; w1:=w4*w3;
Append(~max,sub<G|w1,w2>); w1:=x; w2:=y;
w3:=w1*w2; w4:=w3*w2; w5:=w3*w4; w6:=w3*w5; w8:=w6*w5; w9:=w3*w8; w10:=w9*w4;
w4:=w3^7; w5:=w4*w10; w6:=w10*w4; w7:=w6^-1; w6:=w7*w5; w7:=w6^10; w12:=w10*w7;
w4:=w1*w1; w11:=w4*w3; w3:=w11*w12; w4:=w3*w12; w5:=w3^5; w6:=w5*w5; w3:=w4^-1;
w7:=w4*w6; w6:=w7*w3; w7:=w6^-1; w8:=w7*w5; w5:=w8*w6; w3:=w2^-1; w6:=w5^4;
w7:=w3*w6; w6:=w7*w2; w7:=w4*w5; w8:=w7*w5; w9:=w7*w8; w10:=w7*w9; w8:=w10*w9;
w7:=w8*w8; w8:=w6*w7; w10:=w8^11; w9:=w2*w10; w8:=w9^-1; w10:=w8*w11; w11:=w10*w9;
w10:=w8*w12; w12:=w10*w9; w3:=w11*w12; w6:=w3^8; w7:=w3*w12; w8:=w3*w7;
w7:=w3*w8; w8:=w7^5; w7:=w8*w6; w8:=w7^-1; w9:=w8*w6; w10:=w9*w7;
w1:=w4*w10; w6:=w7*w5; w2:=w8*w6;
Append(~max,sub<G|w1,w2>); w1:=x; w2:=y;
w3:=w1*w2; w4:=w3*w2; w5:=w3*w4; w6:=w3*w5; w7:=w6*w3; w8:=w7*w4;
w2:=w8^3; w3:=w5*w7; w5:=w3^31; w3:=w5^-1; w8:=w3*w1; w1:=w8*w5;
w3:=w4*w6; w4:=w3^13; w3:=w4^-1; w8:=w4*w2; w2:=w8*w3;
Append(~max,sub<G|w1,w2>); w1:=x; w2:=y;
w3:=w1*w2; w4:=w3*w2; w5:=w4*w2; w6:=w5*w2; w7:=w6*w5; w6:=w7*w3; w5:=w4^4;
w7:=w5*w6; w6:=w3*w7; w7:=w6^12; w3:=w2^-1; w4:=w7*w3; w3:=w7*w2; w5:=w4*w3;
w4:=w5^3; w2:=w1*w4; w1:=w7*w6; w7:=w6*w4; w5:=w6*w7; w3:=w2*w5; w7:=w6*w5;
w2:=w3*w7; w3:=w2*w4; w5:=w7*w7; w2:=w3*w5; w7:=w6^5; w3:=w2*w7; w2:=w3*w5;
Append(~max,sub<G|w1,w2>); w1:=x; w2:=y;
w3:=w1*w2; w4:=w3*w2; w5:=w3*w4; w6:=w3*w5; w3:=w4*w2; w2:=w5*w3; w5:=w2*w3;
w2:=w6^17; w3:=w4^21; w7:=w2*w3; w2:=w7^-1; w3:=w2*w1; w1:=w3*w7; w2:=w4^16;
w3:=w6^30; w7:=w2*w3; w6:=w7^-1; w3:=w6*w5; w2:=w3*w7;
Append(~max,sub<G|w1,w2>); w1:=x; w2:=y;
w3:=w1*w2; w4:=w3*w2; w5:=w3*w4; w6:=w3*w5; w7:=w6*w3; w8:=w7^3; w7:=w6^7;
w6:=w7^-1; w4:=w7*w1; w7:=w4*w6; w4:=w5*w3; w5:=w4^15; w4:=w5^-1; w9:=w5*w8;
w8:=w9*w4; w2:=w7*w8; w3:=w2^3; w4:=w1*w3; w1:=w4^20; w3:=w1*w2; w4:=w3*w2;
w5:=w3*w4; w6:=w3*w5; w8:=w6*w5; w9:=w3*w8; w1:=w9^3; w3:=w1*w2; w4:=w3*w2;
w5:=w3*w4; w6:=w3*w5; w9:=w3^3; w10:=w4^5; w8:=w9*w10; w9:=w8^-1; w10:=w9*w2;
w1:=w10*w8; w2:=w6^8; w4:=w3^-1; w5:=w3*w2; w2:=w5*w4; w3:=w1*w2; w4:=w3*w2;
w5:=w3*w4; w6:=w3*w5; w9:=w6*w3; w4:=w3*w3; w3:=w4^-1; w10:=w3*w9;
w2:=w10*w4; w3:=w7*w2; w1:=w3^9;
Append(~max,sub<G|w1,w2>);
\end{verbatim}
\subsection{$Th$}

Here we assume that $G=Th$ is given by standard generators $x$ in class 2A and $y$ in class 3A. All of these subgroups are given in the online ATLAS. This code produces a set \texttt{max} of maximal subgroups of $G$.

This code takes seconds to run on the author's computer.

\begin{verbatim}
OrderG:=90745943887872000;
n2:=OrderG/92897280;
n3:=OrderG/12737088+OrderG/472392+OrderG/174960;
n5:=OrderG/3000;
n7:=OrderG/1176;
n13:=OrderG/39;
n19:=OrderG/19;
n31:=2*OrderG/31;

max:=[]; w1:=x; w2:=y;
w3:=w1*w2; w4:=w3*w2; w5:=w3*w4; w6:=w3*w5; w7:=w6*w3; w8:=w7*w4; w9:=w3*w8;
w8:=w7*w9; w2:=w8^5; w4:=w3^8; w3:=w4^-1; w5:=w3*w1; w1:=w5*w4;
Append(~max,sub<G|w1,w2>); w1:=x; w2:=y;
w3:=w1*w2; w4:=w3*w2; w5:=w4*w3; w6:=w5*w5; w5:=w3^15; w7:=w4^9;
w8:=w5*w7; w5:=w3^12; w7:=w8*w5; w5:=w4^16; w8:=w7*w5; w5:=w3^17;
w7:=w8*w5; w8:=w7^-1; w3:=w8*w6; w2:=w3*w7;
Append(~max,sub<G|w1,w2>); w1:=x; w2:=y;
w3:=w1*w2; w4:=w3*w2; w5:=w3*w4; w6:=w3*w5; w7:=w1*w6; w8:=w6*w1;
w9:=w8^-1; w8:=w9*w7; w9:=w8^9; w2:=w7*w9; w8:=w6*w5; w9:=w3*w8;
w10:=w9* w4; w7:=w1*w10; w8:=w10*w1; w9:=w8^-1; w8:=w9*w7; w1:=w8^14;
Append(~max,sub<G|w1,w2>); w1:=x; w2:=y;
w3:=w1*w2; w4:=w3*w2; w5:=w4^-1; w3:=w5*w2; w2:=w3*w4;
Append(~max,sub<G|w1,w2>); w1:=x; w2:=y;
w3:=w1*w2; w4:=w3*w2; w5:=w3*w4; w6:=w3*w5; w7:=w6*w3; w8:=w7*w4; w9:=w3*w8;
w10:=w9*w4; w8:=w3*w7; w2:=w8^5; w5:=w3*w10; w6:=w5*w4; w7:=w9*w6; w8:=w6*w9;
w6:=w7^27; w7:=w8^7; w8:=w6*w7; w7:=w8^-1; w5:=w7*w1; w1:=w5*w8; w5:=w4*w9;
w6:=w5*w3; w7:=w10*w6; w8:=w7^10; w7:=w8^-1; w5:=w7*w2; w2:=w5*w8;
Append(~max,sub<G|w1,w2>); w1:=x; w2:=y;
w3:=w1*w2; w4:=w3*w2; w5:=w3*w4; w6:=w5*w5; w5:=w6*w4; w6:=w3*w3; w7:=w6*w5;
w2:=w7^3; w5:=w3^8; w6:=w4^6; w7:=w5*w6; w6:=w7^-1; w5:=w6*w1; w1:=w5*w7;
w5:=w4^13; w6:=w3^5; w7:=w5*w6; w6:=w7^-1; w5:=w6*w2; w2:=w5*w7;
Append(~max,sub<G|w1,w2>); w1:=x; w2:=y;
w3:=w1*w2; w4:=w3*w2; w5:=w3^3; w6:=w5*w2; w2:=w6^3; w6:=w4^3;
w4:=w3^11; w3:=w4*w5; w5:=w6*w4; w4:=w6*w6; w6:=w3*w4; w3:=w5^-1;
w4:=w3*w1; w1:=w4*w5; w3:=w6^-1; w4:=w3*w2; w2:=w4*w6;
Append(~max,sub<G|w1,w2>); w1:=x; w2:=y;
w3:=w1*w2; w4:=w3*w2; w5:=w3^4; w3:=w5*w2; w2:=w3*w4; w3:=w4*w4;
w6:=w3*w3; w4:=w5*w3; w3:=w6*w4; w4:=w6^3; w6:=w4*w5; w5:=w3*w4;
w3:=w5^-1; w4:=w3*w1; w1:=w4*w5; w3:=w6^-1; w4:=w3*w2; w2:=w4*w6;
Append(~max,sub<G|w1,w2>); w1:=x; w2:=y;
w3:=w1*w2; w4:=w3*w2; w5:=w3*w4; w6:=w5*w5; w5:=w6*w4; w6:=w3^4;
w2:=w6*w5; w1:=w2^10; w5:=w6*w4; w2:=w5^5; w7:=w3^7; w3:=w6*w6;
w5:=w4^5; w6:=w5*w7; w5:=w3*w6; w6:=w5^-1; w7:=w6*w1; w1:=w7*w5;
w7:=w3*w3; w6:=w4^11; w5:=w7*w6; w6:=w5^-1; w7:=w6*w2; w2:=w7*w5;
Append(~max,sub<G|w1,w2>); w1:=x; w2:=y;
w3:=w1*w2; w4:=w3*w2; w5:=w3^5; w6:=w5*w2; w2:=w6^5; w5:=w3^3;
w6:=w4^15; w7:=w5*w6; w6:=w7^-1; w5:=w6*w1; w1:=w5*w7; w5:=w4^8;
w6:=w3^9; w7:=w5*w6; w6:=w7^-1; w5:=w6*w2; w2:=w5*w7;
Append(~max,sub<G|w1,w2>); w1:=x; w2:=y;
w3:=w1*w2; w4:=w3*w2; w5:=w3^4; w1:=w5*w4; w5:=w4^8; w6:=w5*w4; w4:=w3*w3;
w7:=w5*w4; w4:=w3^5; w3:=w4*w4; w5:=w4*w7; w4:=w3*w6; w3:=w1^5; w6:=w5^-1;
w7:=w2*w5; w1:=w6*w7; w6:=w4^-1; w7:=w3*w4; w2:=w6*w7;
Append(~max,sub<G|w1,w2>); w1:=x; w2:=y;
w3:=w1*w2; w4:=w3*w2; w5:=w3*w4; w6:=w3*w5; w2:=w6^6; w5:=w4^6; w6:=w5^-1;
w4:=w6*w2; w2:=w4*w5; w4:=w3*w3; w3:=w4^-1; w5:=w3*w1; w1:=w5*w4;
Append(~max,sub<G|w1,w2>); w1:=x; w2:=y;
w3:=w1*w2; w4:=w3*w2; w5:=w3*w4; w6:=w3*w5; w2:=w6^6; w5:=w4^8; w4:=w5^-1;
w6:=w4*w1; w1:=w6*w5; w4:=w3^8; w3:=w4^-1; w6:=w3*w2; w2:=w6*w4;
Append(~max,sub<G|w1,w2>); w1:=x; w2:=y;
w3:=w1*w2; w4:=w3*w2; w5:=w3^5; w6:=w5*w2; w2:=w6^5; w5:=w3^4;
w6:=w4^4; w7:=w5*w6; w6:=w7^-1; w5:=w6*w1; w1:=w5*w7; w5:=w4^17;
w6:=w3^6; w7:=w5*w6; w6:=w7^-1; w5:=w6*w2; w2:=w5*w7;
Append(~max,sub<G|w1,w2>); w1:=x; w2:=y;
w3:=w1*w2; w4:=w3*w2; w5:=w2*w3; w11:=w1*w5; w12:=w5*w1; w13:=w12^-1; w12:=w13*w11;
w11:=w12^4; w6:=w5*w11; w11:=w4*w3; w5:=w3*w4; w12:=w11^5; w13:=w5^5; w5:=w12*w13;
w1:=w3*w3; w2:=w1*w1; w7:=w1*w2; w11:=w4*w4; w12:=w11*w11; w13:=w2*w11; w14:=w3*w13;
w11:=w3*w12; w3:=w2*w7; w2:=w1*w3; w1:=w14*w11; w14:=w13*w4; w13:=w12*w14; w11:=w2*w13;
w12:=w4^10; w13:=w3*w12; w4:=w13*w7; w3:=w11^-1; w2:=w5^3; w7:=w2*w6; w11:=w6^4;
w12:=w11*w2; w13:=w2*w11; w14:=w5*w13; w13:=w11*w11; w11:=w5*w6; w2:=w11*w13;
w13:=w5*w14; w5:=w7^4; w6:=w11^6; w11:=w13*w12; w12:=w2*w14; w13:=w12*w11; w14:=w6^-1;
w11:=w5*w14; w12:=w5*w6; w2:=w11*w12; w7:=w2*w13; w11:=w3*w7; w2:=w11*w4;
Append(~max,sub<G|w1,w2>); w1:=x; w2:=y;
w3:=w1*w2; w4:=w3*w2; w5:=w3*w4; w6:=w3*w5; w7:=w6*w3; w8:=w3*w7; w2:=w8^5;
w8:=w7*w4; w9:=w3*w8; w10:=w9*w4; w11:=w4*w9; w12:=w11*w3; w11:=w12*w10; w13:=w10*w12;
w8:=w13^29; w6:=w11^23; w7:=w8*w6; w6:=w7^-1; w8:=w6*w2; w2:=w8*w7; w5:=w9*w3;
w6:=w10*w4; w4:=w5*w6; w3:=w4^27; w4:=w3^-1; w5:=w4*w1; w1:=w5*w3;
Append(~max,sub<G|w1,w2>);
\end{verbatim}
\bibliography{references}
\end{document}